\documentclass{amsart}

\usepackage{amsmath}
\usepackage{amsfonts}
\usepackage{amssymb}
\usepackage{amsthm}
\usepackage{array}
\usepackage{bbm}
\usepackage{chngpage}
\usepackage{comment}
\usepackage{float}
\usepackage[OT2,T1]{fontenc}
\usepackage{graphicx,caption}
\usepackage[export]{adjustbox}
\usepackage{longtable}
\usepackage{mathrsfs}
\usepackage{mathtools}
\usepackage{wrapfig}
\usepackage{cutwin}
\usepackage{shapepar}
\usepackage{tikz}
\usepackage[lowtilde]{url}
\usepackage{subcaption}
\usetikzlibrary{matrix,arrows}
\usepackage{tabularx}
\usepackage{multirow}
\usepackage[colorlinks, citecolor = blue]{hyperref}
\usepackage[capitalize, nameinlink]{cleveref}
\usepackage{enumitem}
\usepackage{diagbox}
\usepackage{etoolbox}
\usepackage{algpseudocode}
\usepackage{etoolbox}

\usepackage[margin=0.9in]{geometry}

\patchcmd{\subsection}{-.5em}{.5em}{}{}

\newcommand{\F}{\mathbb{F}}

\newcommand{\Q}{\mathbb{Q}}
\newcommand{\R}{\mathbb{R}}
\newcommand{\Z}{\mathbb{Z}}

\newcommand{\cF}{\mathcal{F}}

\DeclareSymbolFont{cyrletters}{OT2}{wncyr}{m}{n}
\DeclareMathSymbol{\sha}{\mathalpha}{cyrletters}{"58}

\newlength{\strutheight}
\newcommand{\half}{\frac{1}{2}}
\newcommand{\thalf}{\tfrac{1}{2}}
\newcommand{\tth}{\text{th}}

\newtheorem{theorem}{Theorem}[section]
\newtheorem{lemma}[theorem]{Lemma}

\newtheorem{heuristic}[theorem]{Heuristic}

\theoremstyle{definition}

\newtheorem{figurecap}[theorem]{Figure}
\newtheorem{example}[theorem]{Example}
\AtBeginEnvironment{example}{%
  \pushQED{\qed}%
}
\AtEndEnvironment{example}{\popQED\endexample}

\author{Alex Cowan}
\address{Department of Mathematics, Harvard University, Cambridge, MA 02138 USA}
\email{cowan@math.harvard.edu}
\thanks{The author was supported by the Simons Foundation Collaboration Grant 550031.}

\title{Murmurations and explicit formulas}
\date{\today}

\begin{document}
\begin{abstract}
  Unexpected oscillations in $a_p$ values in a family of elliptic curves were observed experimentally in \cite{HLOP}. We propose a heuristic explanation for these oscillations based on the ``explicit formula'' from analytic number theory. A crucial ingredient in this heuristic is that the distribution of the zeros of the associated $L$-functions has a quasi-periodic structure. We present empirical results for a family of elliptic curves, a family of quadratic Dirichlet characters whose values exhibit similar oscillations, and a family of Dirichlet characters whose values do not.
\end{abstract}
\maketitle
\tableofcontents

\section{Introduction}

\textit{Murmurations} are oscillations in the average $a_p$ values of a set of elliptic curves as $p$ varies. This phenomenon was observed empirically in \cite{HLOP} using techniques from data science. The oscillations are unexpected and quite striking, and it is not immediately clear what causes them. \textit{Explicit formulas} are identities arising in analytic number theory which relate partial sums of the prime-power coefficients of $L$-functions to the zeros of those $L$-functions.

In this paper, we present experimental evidence which suggests that murmurations can be viewed as a natural consequence of averaging an explicit formula. Our main empirical observation is visualized in \cref{fig:E_murmurations}. The connection we outline between murmurations and explicit formulas is purely heuristic: when we average the explicit formula over sets of elliptic curves, we ignore a mysterious error term.

From the perspective of explicit formulas, a critical ingredient for the existence of murmurations is a quasi-periodic structure in the density of zeros in the associated family of $L$-functions. In \cref{fig:zero_dist} we illustrate how two families similar to those considered in \cite{HLOP} exhibit this quasi-periodic structure.

Quasi-periodic structure in $L$-function zeros is expected to exist for many kinds of arithmetic objects, not just elliptic curves. In \cref{section:dirichlet} we present empirical evidence that quadratic Dirichlet characters exhibit murmurations, and those murmurations are similarly accounted for by an explicit formula. In that section, we also present a family of Dirichlet characters whose $L$-function zeros have no clear structure, and which do not exhibit murmurations.

\section{The explicit formula}

Let $E/\Q$ be an elliptic curve with rank $r(E)$, conductor $N(E)$, and $L$-function $L_E$. Let $a_p(E)$ denote the $p^\tth$ Fourier coefficient of the rational newform associated to $E$. For $p \nmid N(E)$ we have $a_p(E) = p+1 - \#E(\F_p)$. Assume the Riemann Hypothesis for $L_E$. For positive integers $n$, define $\gamma_n(E)$ so that $\{\thalf + i\gamma_n(E) \,:\, n \in \Z_{>0}\}$ runs over the nontrivial zeros of $L_E$ with strictly positive imaginary part, and so that $\gamma_n(E) < \gamma_m(E)$ if and only if $n < m$. Define $\gamma_{-n}(E) \coloneqq -\gamma_n(E)$. Our assumptions imply that $\{\thalf + \gamma_n(E) \,:\, n \in \Z_{\neq 0}\}$ is the set of all nontrivial zeros of $L_E(s)$ with nonzero imaginary part.

By evaluating the inverse Mellin transform
\begin{align*}
  \int_{\sigma - i\infty}^{\sigma + i\infty} \frac{L_E'(s)}{L_E(s)}X^s\frac{ds}{s}
\end{align*}
in two different ways, first by writing the logarithmic derivative $L_E'(s)/L_E(s)$ as a Dirichlet series, and second in terms of the residues of the integrand, one obtains the following \textit{explicit formula}.
\begin{lemma}[{\cite[Lemma 2.2]{fiorilli}}]\label{explicit_formula}
  \begin{align*}
    \frac{\log x}{x^\half}\sum_{p<x} \frac{a_p(E)}{p^\half} = 1 - 2r(E) - \sum_{n \neq 0} \frac{x^{i\gamma_n(E)}}{\half + i\gamma_n(E)} + \mathrm{Err}_E(x)
  \end{align*}
  with $\mathrm{Err}_E(x) = o_E(1)$.
\end{lemma}

A thorough treatment of the explicit formula for the Riemann zeta function and for primitive Dirichlet $L$-functions can be found in \cite[Ch.\ 12]{MV}. Treatment yielding more precise descriptions of the error term $\mathrm{Err}_E(x)$ on the right hand side above can be found in \cite[Lemma 2.2]{fiorilli} as well as \cite[(2.7)]{kim_murty}, \cite[\S 6]{conrad}, \cite[\S 4]{devin}, etc. In particular, \cite[(1.11)]{rubinstein:2013} highlights that the error term coming from the handling of the Dirichlet series is a sum over prime powers, and particularly squares. This is visually apparent in \cref{fig:E_murmurations}.

Let $\cF$ be a finite set of elliptic curve isogeny classes. We will call $\cF$ a \textit{family}. Let $\delta$ denote the Dirac delta function. Define $\rho_\cF$ to be the distribution
\begin{align*}
  \rho_\cF(\gamma) \coloneqq \frac{1}{\#\cF}\sum_{E \in \cF} \sum_{n \neq 0} \delta(\gamma - \gamma_n(E)).
\end{align*}
Equivalently, $\rho_\cF$ is the distribution satisfying, for every subset $U \subseteq \R$,
\begin{align*}
  \rho_\cF(U) = \frac{1}{\#\cF}\cdot\#\big\{(E,n) \in \cF \times \Z_{\neq 0} \,:\, \gamma_n(E) \in U \big\}.
\end{align*}
\cref{fig:zero_dist} visualizes $\rho_\cF$, the distribution of $L$-function zeros of curves in $\cF$, for two families $\cF$ similar to the ones considered in \cite{HLOP}. For certain families $\cF$, e.g.\ large and containing similar curves, it may be reasonable to approximate $\rho_\cF$ by a density, i.e.\ a real-valued function, in the sense of weak convergence.

We propose the following heuristic, which results from averaging the explicit formula given in \cref{explicit_formula} over $\cF$ and neglecting error terms.
\begin{heuristic}\label{avg_explicit_formula}
  \begin{align*}
    \frac{\log x}{x^\half}\sum_{p<x} \frac{1}{\#\cF}\sum_{E\in\cF} \frac{a_p(E)}{p^\half} \approx 1 - \frac{2}{\#\cF}\sum_{E\in\cF}r(E) - \int_\R \frac{\rho_\cF(\gamma)}{\half + i\gamma} e^{i\gamma \log x}\,d\gamma.
  \end{align*}
\end{heuristic}

The integral on the right hand side above can be thought of as an inverse Fourier transform of the distribution $\rho_\cF(\gamma)(\thalf + i\gamma)^{-1}$ with space variable $\log x$.

The error term which is neglected in passing from \cref{explicit_formula} to \cref{avg_explicit_formula} is, in the notation of \ref{explicit_formula},
\begin{align}\label{avg_error}
  \frac{1}{\#\cF}\sum_{E \in \cF} \mathrm{Err}_E(x).
\end{align}
The discussion following \cref{explicit_formula} above lists a couple references which describe $\mathrm{Err}_E(x)$ in various ways. However, at the moment we find the average in \eqref{avg_error} mysterious. 
It is conceivable to us that one could produce reasonable heuristic bounds for the average \eqref{avg_error}, in particular by accounting for the prime powers which arise naturally in explicit formulas, but we won't pursue this here.%
\section{Empirical results}\label{section_empirical}

\cref{fig:zero_dist} illustrates the distribution $\#\cF \cdot \rho_\cF$ of low-lying zeros for two choices of families $\cF$ which are similar to the ones considered in \cite{HLOP}.
\begin{figure}[H]
  \includegraphics[width=\textwidth]{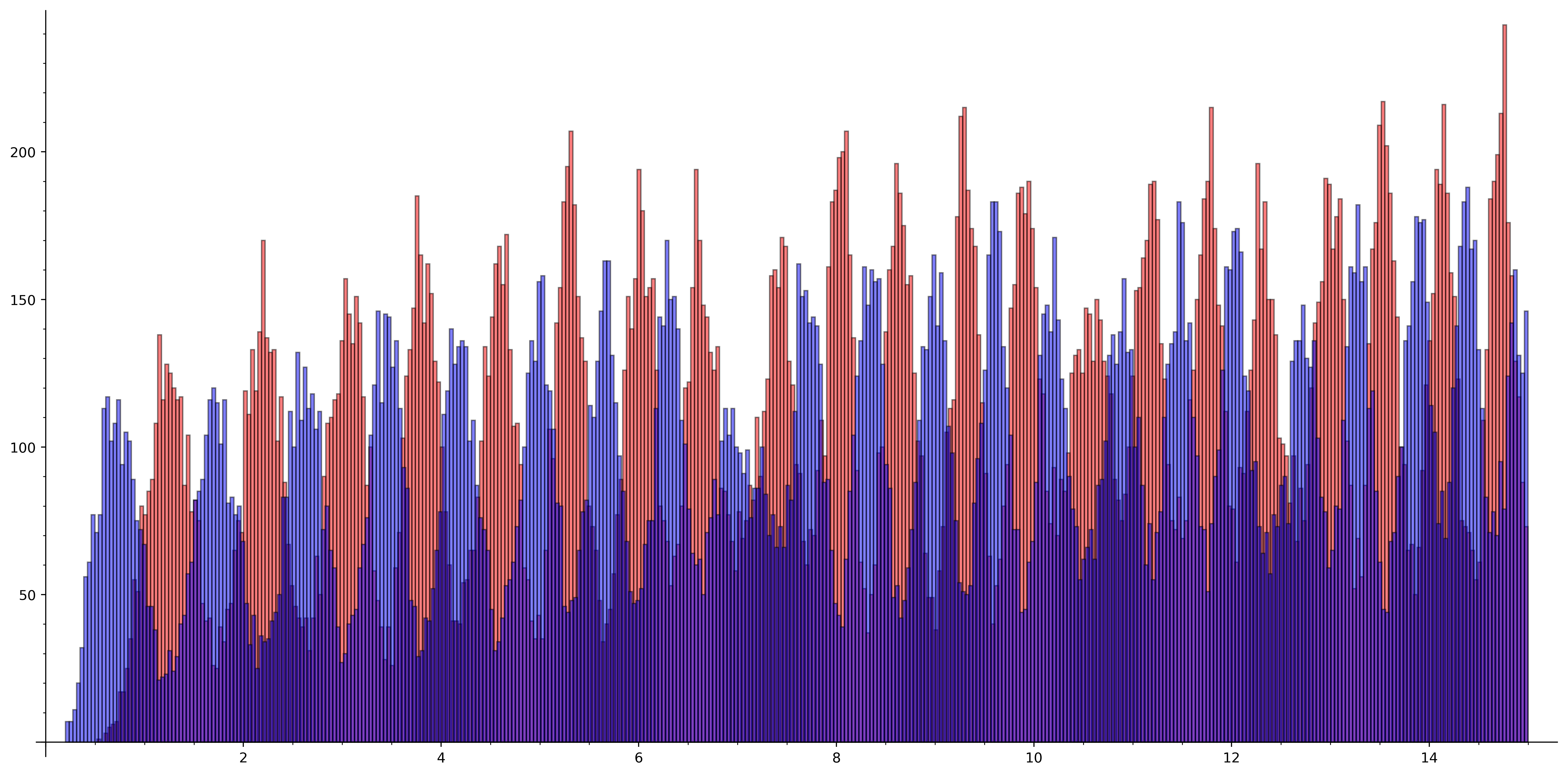}
\end{figure}
\begin{figurecap}\label{fig:zero_dist}
  Histogram of $0 < \gamma_n(E) < 15$ for all elliptic curve isogeny classes with conductor in $[9000, 10000]$ and rank $0$ (blue; 1784 curves) or $1$ (red; 2118 curves).
\end{figurecap}

The visually obvious structure in both distributions $\rho_\cF$ depicted in \cref{fig:zero_dist} leads to the inverse Fourier transforms of $\rho_\cF(\gamma)(\thalf + i\gamma)^{-1}$ looking periodic-ish, in a way which is hard to formulate rigorously but striking in the same way the data from \cite{HLOP} is. For the set of elliptic curve isogeny classes with rank $0$ and conductor between $9000$ and $10000$ (the blue family in \cref{fig:zero_dist}), \cref{fig:E_murmurations} shows this inverse Fourier transform, the left hand side of \cref{avg_explicit_formula}, and their sum.%
\begin{figure}[H]
  \includegraphics[width=0.9\textwidth]{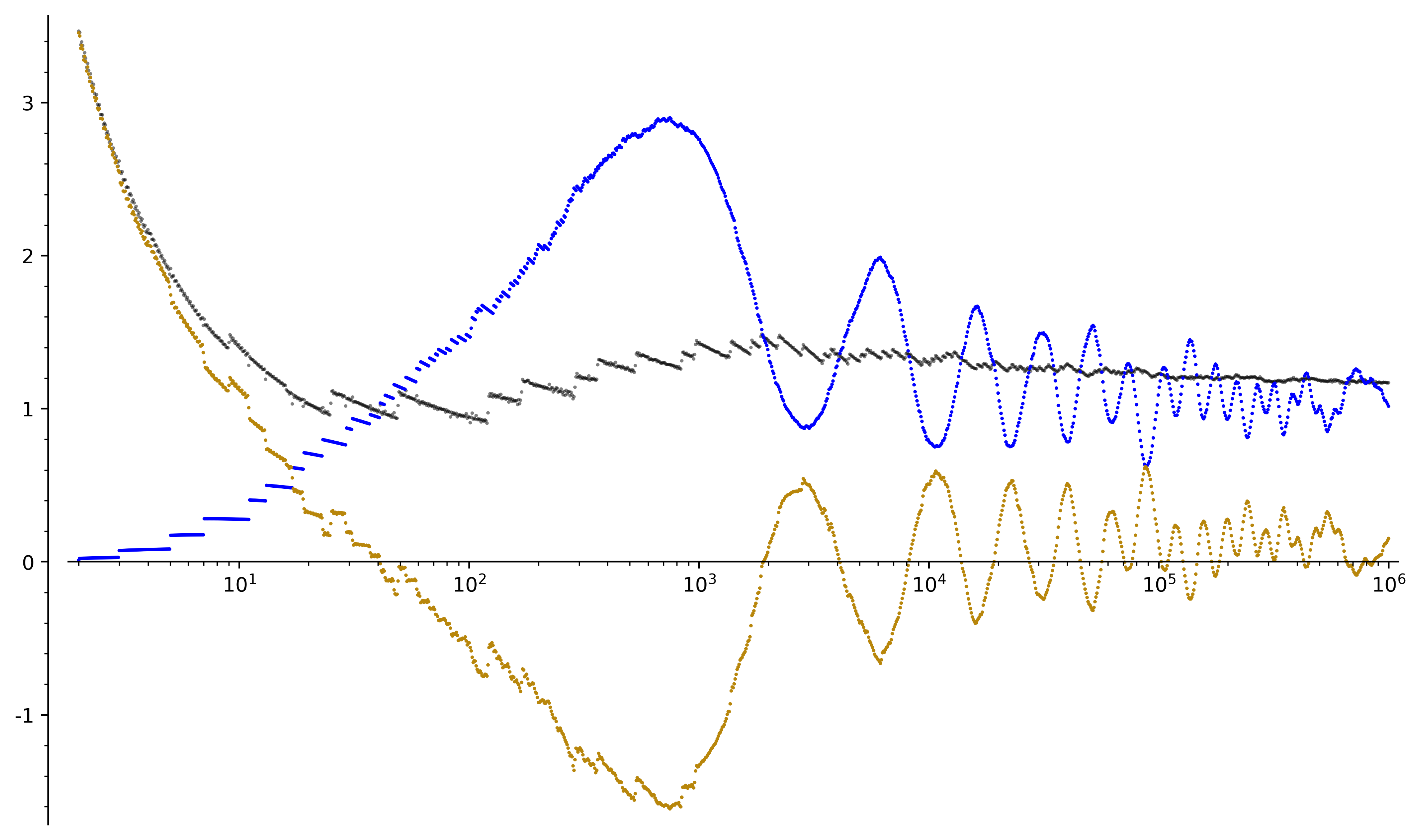}
\end{figure}
\begin{figurecap}\label{fig:E_murmurations}
  For the set of elliptic curve isogeny classes with conductor in $[9000, 10000]$ and rank $0$, left hand side of \ref{avg_explicit_formula} (blue), integral on the right hand side of \cref{avg_explicit_formula} (gold) using only the $500$ lowest-lying zeros with positive imaginary part and their negatives for each curve, and the sum of the blue and gold values (black).
\end{figurecap}

In the figure above, the black curve, which is $1$ plus the average error \eqref{avg_error}, shows clear jumps at squares of primes and at $16$. This is explained by the fact that explicit formulas sum over all prime powers, see \cite[(2.4)]{kim_murty} or the other references listed in the discussion following \cref{explicit_formula}. The black curve in \cref{fig:E_murmurations} looks to fluctuate very little otherwise, suggesting to us that \cref{avg_explicit_formula} would be quite reasonable if modified to account for the aforementioned prime powers.

The handful of isolated black data points lying just off the curve arise when evaluating \cref{avg_explicit_formula} for $x$ nearly a prime. The left hand side is discontinuous at prime $x$, while the right hand side, if the sum over zeros is truncated, is continuous for all $x$. A version of Perron's formula that bounds the error term when $x$ is nearly integral is given in \cite[(2.8)]{kim_murty} and many other places, and the error term induced by truncating the zero sum is given in \cite[Lemma 2.2]{fiorilli}.%
\section{Murmurations for Kronecker symbols}\label{section:dirichlet}

Explicit formulas exist for $L$-functions in general, see \cite[\S 5.5]{IK}. In particular, the explicit formula for Dirichlet characters is given in \cite[Cor.\ 12.11]{MV}. We will use it in the following form.
\begin{lemma}\label{explicit_formula_dirichlet}
  For any even primitive nontrivial Dirichlet character $\chi$ mod $N$,
  \begin{align*}
    \frac{1}{x^\half}\sum_{p<x} \chi(p)\log p = -\frac{\log x}{x^\half} - \sum_\gamma \frac{x^{i\gamma}}{\half + i\gamma} + R_\chi(x),
  \end{align*}
  where the sum is over the imaginary parts of the nontrivial zeros of $L(s,\chi)$,
  \begin{align*}
    R_\chi(x) = \frac{1}{x^\half}\left(\frac{L'(1,\bar\chi)}{L(1,\bar\chi)} + \log\frac{N}{2\pi} - C_0 - \log\sqrt{1 - x^{-2}} - \sum_{k=2}^\infty \sum_{p^k<x} \chi(p^k)\log p\right),
  \end{align*}
  and $C_0 = 0.577\dots$ is the Euler–-Mascheroni constant.
\end{lemma}

Moreover, a wide variety of arithmetic objects are predicted to have low-lying zero structure similar in nature to that observed in \cref{fig:zero_dist}, see e.g.\ \cite{odlyzko, katz_sarnak, rubinstein, ILS, HKS}. In \cref{fig:zero_dist_dirichlet} we show the distribution of low-lying zeros for the $L$-functions attached to the $307$ Kronecker symbols $\left(\frac{D}{\cdot}\right)$ with $D$ a fundamental discriminant between $9000$ and $10000$. In \cref{fig:dirichlet_murmurations}, we show averages of \cref{explicit_formula_dirichlet} in a way which parallels \cref{avg_explicit_formula}.
\begin{figure}[H]
  \includegraphics[width=\textwidth]{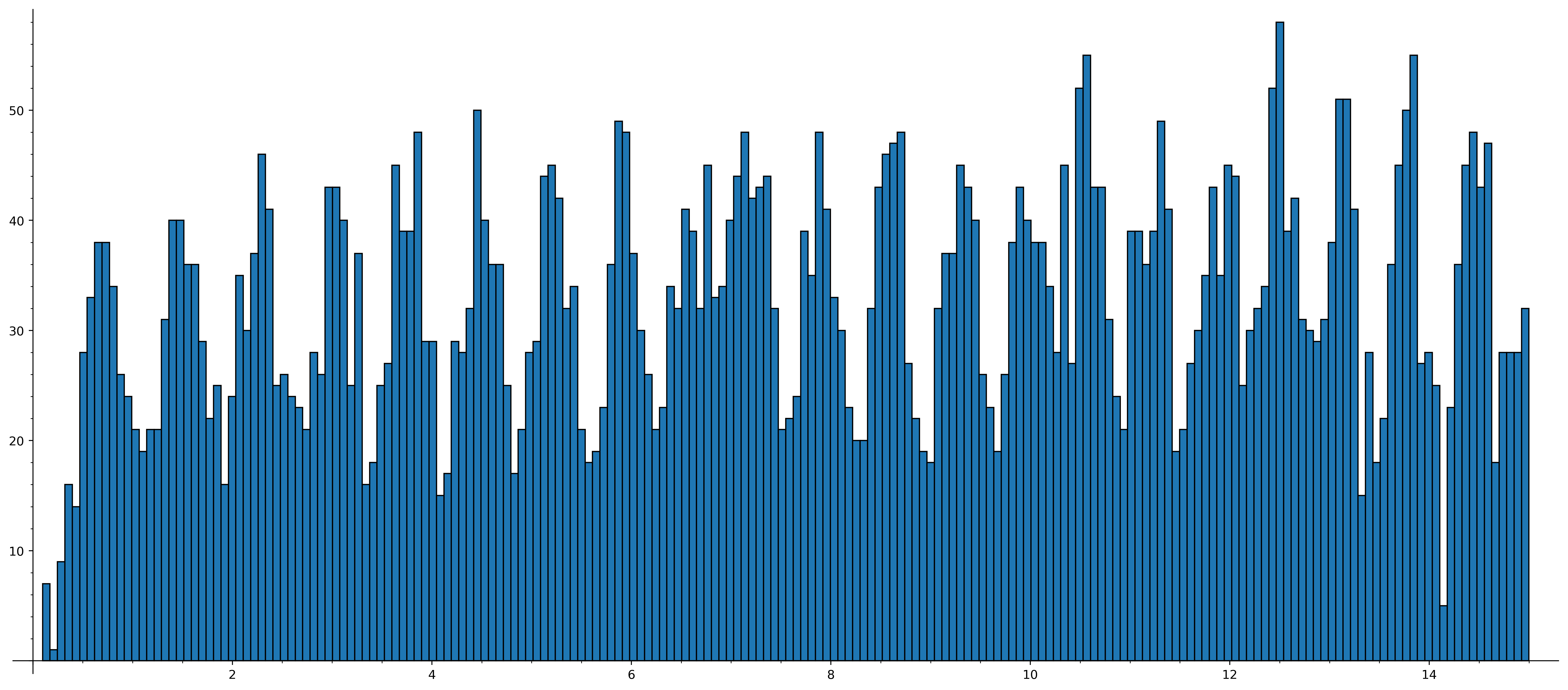}
\end{figure}
\begin{figurecap}\label{fig:zero_dist_dirichlet}
  Histogram of imaginary parts of non-trivial $L$-function zeros for the $307$ Kronecker symbols $\left(\frac{D}{\cdot}\right)$ with $D$ a fundamental discriminant between $9000$ and $10000$.
\end{figurecap}
\begin{figure}[H]
  \includegraphics[width=\textwidth]{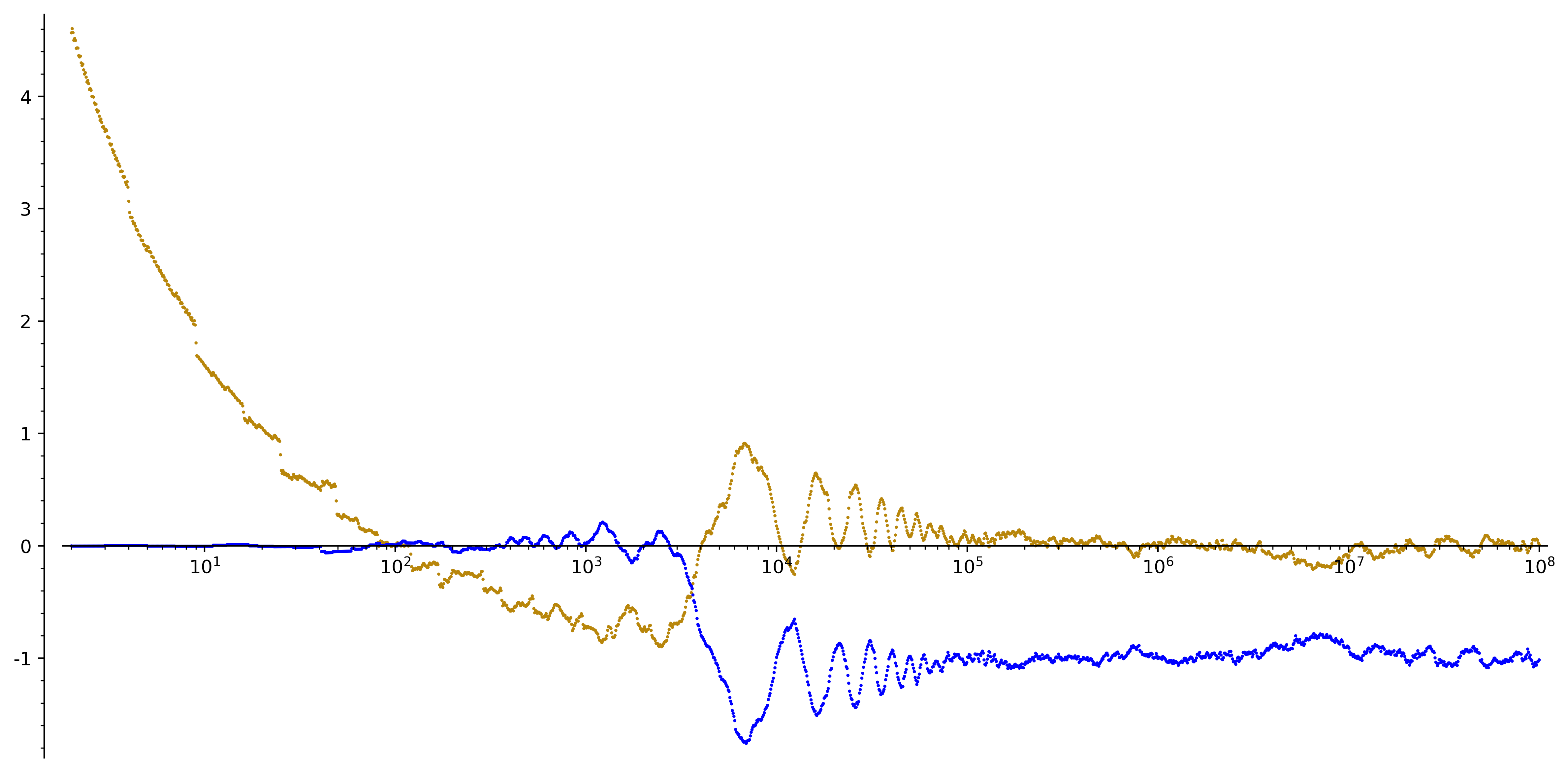}
\end{figure}
\begin{figurecap}\label{fig:dirichlet_murmurations}
  For the $307$ Kronecker symbols $\left(\frac{D}{\cdot}\right)$ with $D$ a fundamental discriminant between $9000$ and $10000$, the left hand side of \ref{explicit_formula_dirichlet} averaged over $D$ (blue), the negative of the right hand side of \cref{explicit_formula_dirichlet} averaged over $D$ (gold) omitting $R_\chi$ and using only zeros with imaginary part at most $200$ in absolute value.
\end{figurecap}

\begin{figure}[H]
\begin{minipage}[t]{0.49\linewidth}
  \centering
  \includegraphics[width=\linewidth, valign=t,left]{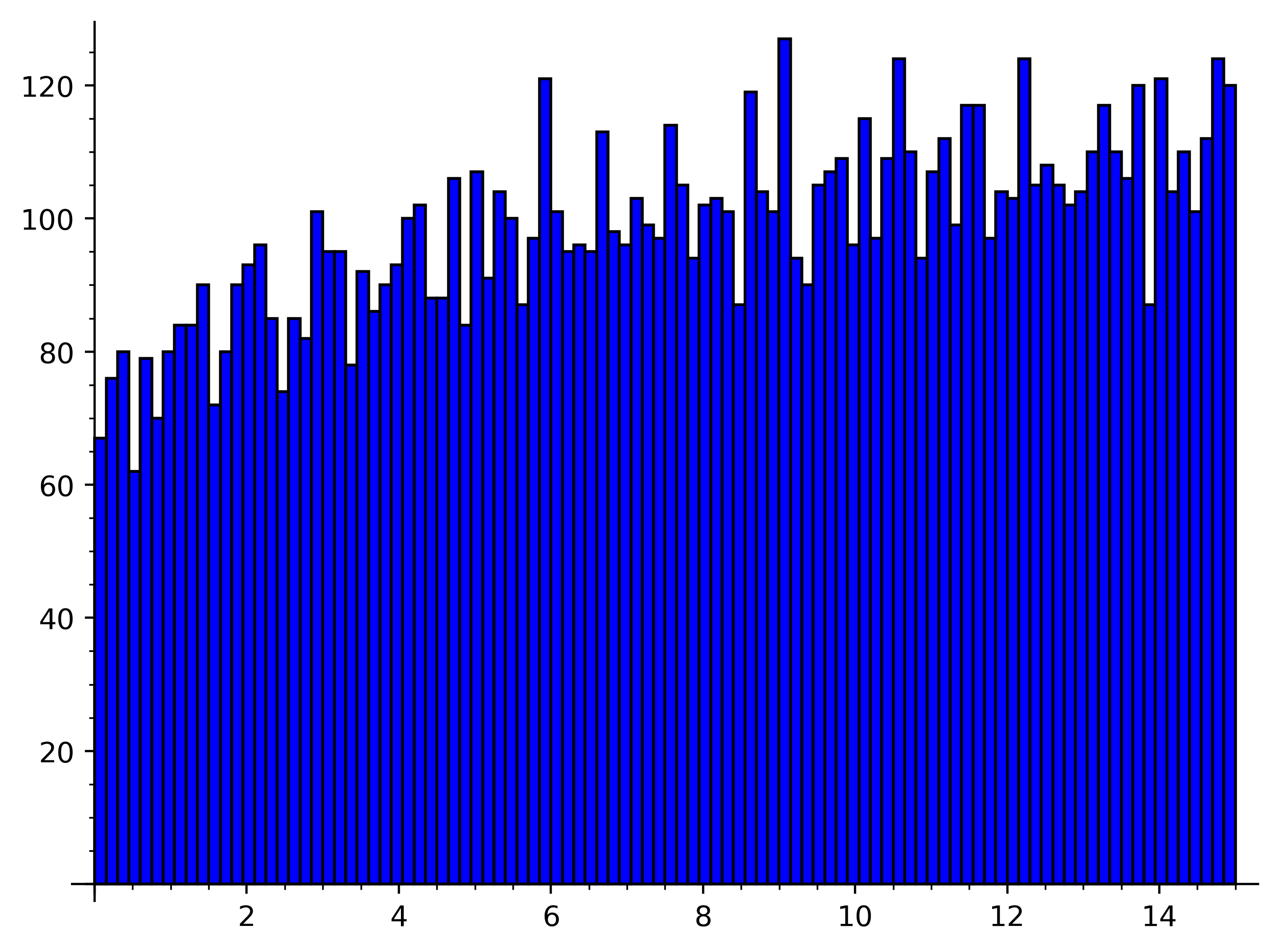}
  \begin{figurecap}\label{fig:zero_dist_2797}
    Histogram of imaginary parts of non-trivial $L$-function zeros for a set of $526$ odd Dirichlet characters mod $2797$ closed under complex conjugation.
  \end{figurecap}
  \vspace{-\baselineskip}
\end{minipage}%
\hfill
\begin{minipage}[t]{0.49\linewidth}
  \centering
  \includegraphics[width=\linewidth, valign=t,right]{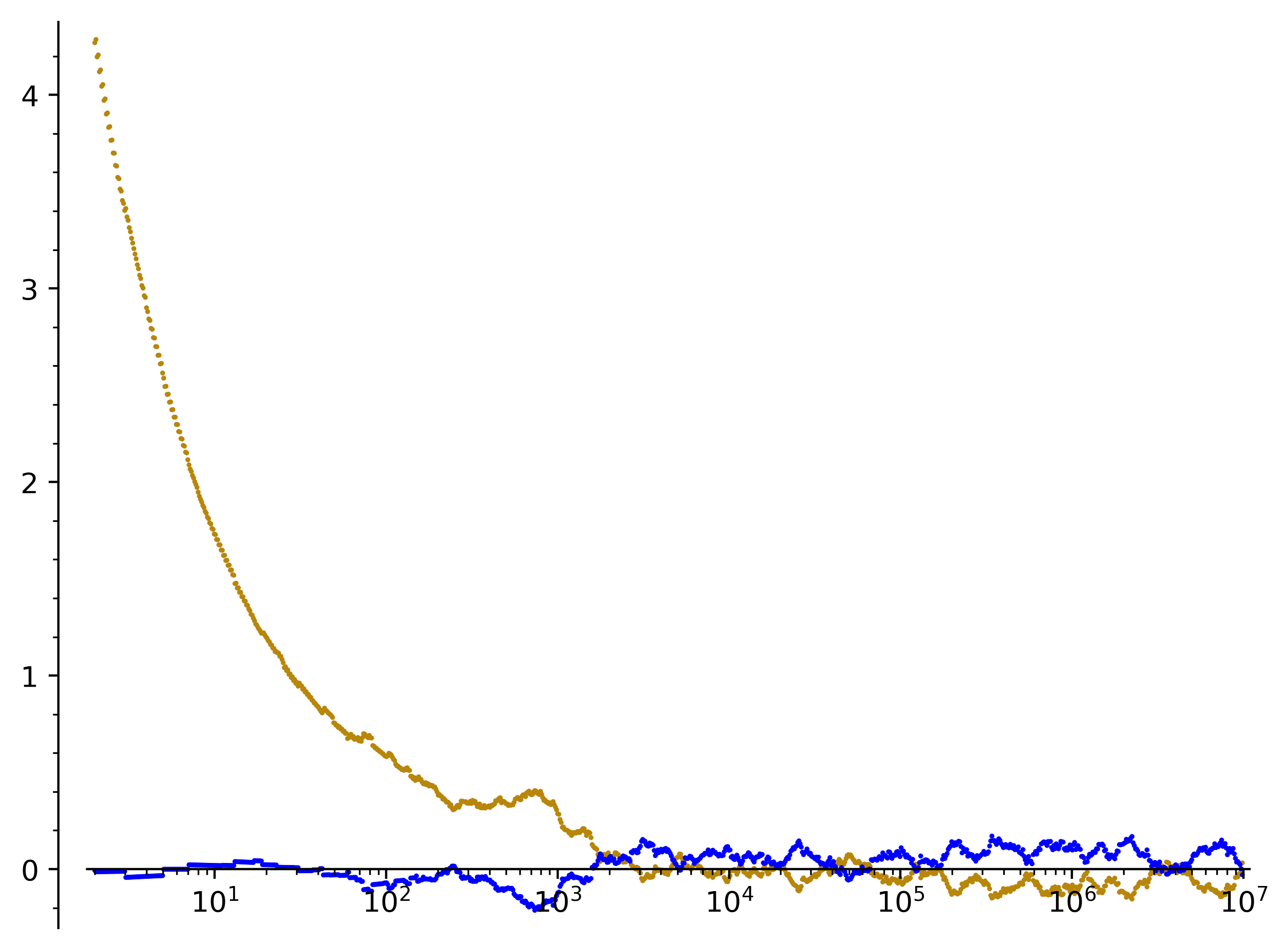}
  \begin{figurecap}\label{fig:murmurations_2797}
    For the Dirichlet characters from \cref{fig:zero_dist_2797}, the left hand side of \ref{explicit_formula_dirichlet} averaged over $D$ (blue), the negative of the sum on the right hand side of \cref{explicit_formula_dirichlet} averaged over $D$ (gold) using only zeros with imaginary part at most $200$ in absolute value.
  \end{figurecap}
  \vspace{-\baselineskip}
\end{minipage}
\end{figure}

In \cref{fig:zero_dist_2797} below we show the low-lying zeros of $L$-functions associated to a randomly chosen set of $526$ odd Dirichlet characters modulo $2797$ closed under complex conjugation. In \cref{fig:murmurations_2797} we show averages of the analogue of \cref{explicit_formula_dirichlet} for odd characters. In contrast to the sets of elliptic curves and Kronecker symbols discussed previously, for these Dirichlet characters there is no clear structure in the low-lying $L$-function zeros, and no clear murmurations. This is consistent with our proposal that the structure in the locations of the zeros is what causes murmurations.


\section*{Acknowledgements}
We thank Drew Sutherland for his helpful suggestions and feedback.

\bibliographystyle{alpha}
\bibliography{murmurationsbib}{}

\end{document}